\documentclass[12pt,reqno]{amsart}

\newcommand{\ver}{{\it 1}}

\usepackage{latexsym}
\usepackage{graphicx}

\newcommand{\PP}{{\mathbb P}}

\newcommand{\R}{{\mathbb R}}
\newcommand{\C}{{\mathbb C}}

\newcommand{\CP}{\C\PP}

\newcommand{\E}{{\mathbf E}}

\renewcommand{\phi}{\varphi}

\newcommand{\dcal}{\mathcal{D}}
\newcommand{\ecal}{\mathcal{E}}

\newcommand{\ical}{\mathcal{I}}

\newcommand{\ocal}{\mathcal{O}}

\newcommand{\rcal}{\mathcal{R}}

\newcommand{\Ga}{\Gamma}

\newcommand{\om}{\omega}
\newcommand{\Om}{\Omega}

\renewcommand{\phi}{\varphi}

\newtheorem{theo}{{\sc Theorem}}[section]

\newtheorem{lem}[theo]{{\sc Lemma}}

\newenvironment{rem}{\medskip\noindent{\it Remark:\/} }{\medskip}

\title[Energy of zeros of random sections on Riemann Surface \ver]{Energy of zeros of random sections on Riemann Surface}

\author{ Qi Zhong}

\address{Department of Mathematics, Johns Hopkins University, USA, 21218}

\email{qzhong@math.jhu.edu}

\thanks{The author would like to thank his advisor Prof. Zelditch for his helpful suggestion on this work. He also wants to
thank Jian Kong and Brain Macdonald for their help on computer
program.}

\begin{document}

\maketitle

\begin{abstract}
The purpose of this paper is to determine the asymptotic of the
average energy of a configuration of $N$ zeros of system of random
polynomials  of degree $N$ as $N \to \infty$ and more generally
the zeros of random holomorphic sections of a line bundle $L\to M$
over any Riemann surface. And we compare our results to the
well-known minimum of energies.
\end{abstract}

\section{Introduction}

 This article is concerned with the asymptotic of the
average energy of the configuration of zeros of $N$-degree random
polynomials as $N \to \infty$ and more generally
 the zeros of random holomorphic sections of a line
bundle $L\to M$ over any compact Riemann surface without boundary.
The energy of a configuration of points $\{z_1, \dots, z_N\}$ on a
surface $M$ equipped with a
 Riemannian  metric $g$ is defined by
\begin{equation}
\ecal^N_{G_g}=\sum_{i\neq j} G_g(z_i, z_j),
 \end{equation} where $G_g$ is the
Green's function for $g$, $G_g(z,w)= -\frac1{2\pi}\chi(z,w)\log
r_g(z,w)+F(z,w)$,where  $F\in C^\infty(M\times M)$ and $\chi(z,w)$
is the cut-off function near the diagonal, we will discuss the
notations in \S 2.5 ; other energies will also be studied.
Electrons moving freely on the surface distribute themselves in a
minimal energy configuration, and many articles have been devoted
to finding the minimal energy configurations and the asymptotic of
the minimal energy.

\bigskip

The question studied in this article is the extent to which zeros
of random polynomials of degree $N$
 tend to resemble minimal energy configurations of N points. Zeros of
random polynomials in complex dimension one repel and like minimal
energy configurations tend to stay $1/\sqrt{N}$ apart. Our main
results show that the average energy of such random zeros is of
the same order of magnitude as that of  minimal energy
configurations.

\bigskip

To state our results, we need some notation. Throughout the
article we identify polynomials of degree $N$ with holomorphic
sections $H^0(\CP^1,\ocal (N))$ of the $N$th power of the
hyperplane section bundle over the complex projective line
$\CP^1$. Our methods apply equally to holomorphic sections
$H^0(M,L^N)$ of powers of a positive holomorphic line bundle $L\to
M$ over any compact Riemann surface.Thus, in addition to studying
zeros of polynomials, we study zeros of random theta functions
over a Riemann surface of genus one, and zeros of random
holomorphic k-differentials over a surface of higher genus.
Moreover, our results apply to general k$\ddot{a}$hler metrics $g$
on these Riemann surfaces.

As recalled in $\S 2.1$, a choice of hermitian metric on L
determines an inner product on $H^0(M, L^N)$ and then a Gaussian
measure $\mu_{h^N}$ on this spaces. Roughly speaking, a random
section $S_c$ is expressed in terms of an orthonormal basis {$S_j$}
of $H^0(M, L^N)$ as $S_c=\sum c_j S_j$ where the $c_j$ are
independent complex normal Gaussian random variables. Also we define
Riemannian form $dV_M=\om_h=\frac i2\Theta_h$. The metric $g$
defining Green's function is not necessarily equal to the metric
derived by $h$.

In the case of $\CP^1$, we also consider the energies $\ecal_s$
defined by
\begin{equation}
\ecal_s(z_1,\dots,z_N)=\sum_{i \not=
j}\frac1{[z_i,z_j]^s}.\end{equation}  In the case $s=0$, the
$s$-energy is defined to be the logarithmic energy
\begin{equation}
\ecal_0(z_1, \dots, z_N) = \sum_{i \not= j} -\log
[z_i,z_j].\end{equation} Here $[z,w]$ is the chordal distance
between two points on $S^2$, where $z$ and $w$ are the two points
on $\CP^1$ corresponding to some points on $S^2$. If $r$ is the
round distance on $S^2$, then the relation between $[\cdot,\cdot]$
and $r$ is
\begin{equation}\label{relat}
[a,b]=\sqrt{2(1-\cos r(a,b))}
\end{equation}
here, $a$, $b$ are points on $S^2$.
\bigskip

We define $\E\ecal^N_{G_g}$ to be the expected (average) value of
the energy of the zeros of Gaussian random sections chosen from the
ensemble $(H^0(M,L^N),\mu_h)$ by
\begin{equation}
\E\ecal^N_{G_g}=\int_{H^0(M,L^N)}\ecal_{G_g}(Z_s)d\mu_h(s)=\int_{\C^{d_N}}\ecal_{G_g}(Z_s)\frac{e^{-|c|^2}}{\pi^{d_N}}dc,
\end{equation}
here $$\ecal_{G_g}(Z_s) = \ecal_{G_g}(z_1, \dots,
z_N)=\sum_{\substack{z_i\neq z_j\\z_i,z_j\in Z_s}}{G_g}(z_i,z_j)$$
where $Z_s = \{z_1, \dots, z_N\}$ is the zeros of $s$ and
$d_N=\dim H^0(M,L^N)$.

Note that if $s$ has a double zero, the energy is infinite, but
this occurs with measure zero.

\bigskip

Recall a Green's function $G_g$ on compact Riemann Manifold
$(M,g)$ without boundary is the kernel of $-\Delta^{-1}_g$. Then
the expected (average) energy is satisfies:
\begin{theo}
 Let $(L,h)\to
M$ be a positive Hermitian line bundle on $M$ with $\frac
i2\Theta_h=\om$, $\Theta$ is the curvature form and $\om$ gives
$M$ a Riemannian volume form. We assume the Chern class of $L$,
$c_1(L)=1$. Then the expected average of Green's function energy
of zeros of random sections is given by:
\begin{equation}\label{main theorem}\E \ecal^N_{G_g}= -\frac1{4\pi}N\log
N+O(N).
\end{equation}
\end{theo}
\begin{rem}
If $c_1(L)\neq 1$, then the number of zeros is $c_1(L)\cdot N.$
\end{rem}

 In [Hr], Elkies proved that
 $$\sum_{i\neq j}g_v(z_j,z_k)\leq \frac12 N\log N+\frac{11}{3}N+o(N).$$
 where $g_v$ is the Green's function with respected to a special
 volume form $d\mu_v$ see [Hr].
 We notice that Elkies' normalization of Green's function is
$$g_v(z,w)=\log {|z,w|_v}+F'(z,w)$$
 which is negative near the diagonal. While, our normalization for
 Green's function
$$G_g(z,w)= -\frac1{2\pi}\chi(z,w)\log
r_g(z,w)+F(z,w)$$
 is positive near the diagonal.
 Then we can rewrite Elkies' result
 \begin{equation}\label{Hr}
(-\frac{1}{2\pi}\sum_{i\neq j}g_v(z_i,z_j))\geq -\frac
1{4\pi}N\log N-\frac {11}{6\pi}N+o(N)
 \end{equation}

\begin{rem}
\begin{enumerate}
\item We see that the leading order term in equation (\ref{Hr}) is
the same as the one in equation (\ref{main theorem}). It means that
the probability that the energy is above the minimum goes to zero as
$N\to\infty$, i.e.
$$P(s:\ecal_{G_g}(s)\leq a+\epsilon N\log N)\leq \frac{o(N)}{\epsilon},$$
where $a=-\frac 1{4\pi}N\log N-\frac {11}{6\pi}N+o(N)$ which is the
minimum of the energy. Above formula is not hard to verify. If
$a=\inf_{s}\ecal_{G_g}(s)$ and $MAX=\sup_{s}\ecal_{G_g}(s)$, then we
get

\begin{align*}
a+o(N)a=& \int_{H^0(M,L^N)}\ecal_{G_g}(Z_s)d\mu_h(s)\notag\\
=& \int^{MAX}_{a}Xd\mu_{\ecal_{G_g}(s)}\notag\\
=&\int^{a+\epsilon a}_a
Xd\mu_{\ecal_{G_g}(s)}+\int^{MAX}_{a+\epsilon
a}Xd\mu_{\ecal_{G_g}(s)}\\
\geq & a\cdot(1-\mu(s:\ecal_{G_g}(s)\geq a+\epsilon a))+(a+\epsilon
a)\cdot\mu(s:\ecal _{G_g}(s)\geq a+\epsilon a)\\
=&a+\epsilon a\cdot\mu(s:\ecal_{G_g}(s)\geq a+\epsilon a)
\end{align*}
So we have $\mu(s:\ecal_{G_g}(s)\geq a+\epsilon a)\leq \frac
{o(N)}{\epsilon}.$

\item Our expect average of Green's function energy is scale
metric invariant,that is, if we rescale the metric $g\to rg$, then
our result (\ref{main theorem}) doesn't change. When $g\to rg$,
$\Delta_g$ operator becomes $\frac 1r\Delta_g$, as we discuss in
\S 2.5, $G_g(z,w)$ is the kernel of  $(-\Delta^{-1})$, then
$G_g(z,w)dV_g\to
 rG_g(z,w)dV_g$ as $g\to rg$. On the other side, $dV_g\to rdV_g$ as
 $g\to rg$, therefore, $G_g(z,w)$ doesn't change as $g\to rg$.

\item The leading term order term is independent of $g$ and $h$.

\item We define the Green's function $G$ to be positive near the
diagonal and the average is negative. Then we conclude that the
off diagonal part dominates the energy.

\end{enumerate}
\end{rem}

\begin{theo}
\begin{enumerate}

\item Consider $\CP^1$ with the Fubini-Study metric $g$, let
$(L,h)\to\CP^1$ be a positive Hermitian line bundle on $\CP^1$
with $\frac i2\Theta_h=\om_g$, $\Theta$ is the curvature form. we
recall equation (2) and have expected  s-energy:
\begin{itemize}
\item when $s=2$

\begin{align}\label{Z2}\E\ecal^N_2=&\frac14{N^2}\log N+\frac {3N^2}4\log(\log N)+\frac
{N^2}2\log 2+\frac12N^2(\frac ML)^2\notag\\&-2N^2\log\frac
ML+o(N^2).
\end{align}

\item When $s<2$
\begin{equation}\label{Z<}\E\ecal^N_{s<2}=\frac{2^{1-s}}{2-s}N^2+\frac
1{(2(2-s))}N^{1+\frac s2}(\log N)^{1-\frac s2}+o(N^{1+\frac
s2}(\log N)^{1-\frac s2})\end{equation}

\item When $2<s<4$ \begin{equation}\label{Z4}
\E\ecal^N_{2<s<4}=C\frac{N^{1+\frac s2}}{4-s}+O(N^{1+\frac s2}(\log
N)^{1-\frac s2}).\end{equation} We will discuss the constant $C$ in
the remark at the end of this section.
\end{itemize}
\item Under the same condition as above,we recall equation (3)
 and have expected logarithmic energy:
\begin{align}
\E \ecal^N_0=&-(\log2-\frac12)N^2+\frac
N2\log^2N-\frac12N\log(\log N)\log N\notag\\
&+\frac12N\log N+\frac12(\log2+1)N +o(N).
\end{align}

\end{enumerate}
\end{theo}

Let us compare our results on average energy to the prior results
on minimal energy. For $s-$energy case, Saff-Kuijilaars in [KS]
identified $\CP^1$ as $S^2\in \R^3$ and
 considered the energy
$$\ecal'_s(x_1, \dots, x_N) = \sum_{i < j} \frac 1{
|x_i-x_j|^s}$$ where $x_i$ are the points on $S^2\subset\R^3$ not
on $\CP^1$ and $|x-y|$ is the chordal distance of $S^2$ . They
investigated the energy
 $\ecal'_s$.
Moreover, they define the minimal $s-$energy for $N$ points on the
sphere
$$E_s(N):=\min_{\{x_1,...,x_N\}}\{\ecal'_s\}.$$

It was proved by Saff-Kuijlaars that when $s=2$, then

\begin{equation}{\label{K2}}E_2 (S^2, N) \sim \frac{1}{8} \; N^2 \log N.\end{equation} And when
$s>2$, then \begin{equation}\label{K>}C_1 N^{1+s/2}\leq E_s(N)\leq
C_2 N^{1+s/2},\end{equation} $C_1,C_2>0$. And when $s<2$, then
\begin{equation}{\label{K<}}
E_2 (S^2, N) \leq \frac12V_2(s)N^2-CN^{1+\frac s2},
\end{equation}
here, $C>0$ and
$V_2(s)=\frac{\Ga(\frac{3}2)\Ga(2-s)}{\Ga((2-s+1)/2)\Ga(2-\frac
s2)}$.

 B.Bergersen, D. Boal and P. Palffy-Muboray in [BBP] identified $\CP^1$ as $S^2\in \R^3$ and
 considered the energy
$$\ecal'_0(x_1, \dots, x_N) = \sum_{i < j} -\log
|x_i-x_j|.$$  They investigated the ground-state energy of
 the logarithm energy of $N$ points $\{x_1,...,x_N\}$, which is the minimal
energy of $\ecal'_0$ for large $N$:
\begin{equation}\label{BBP2}
\min_{\{x_1,...,x_N\}}\{\ecal'_0\}=-(\frac{1}{2}\log 2-\frac
1{4})N^2-\frac{N}{4}\log N+\cdots^1
\end{equation}
\footnotetext[1]{We use $E_o$ to consist our notation, in [BBP],
they defined their own notation.} And in that paper, they gave a
formula for the ground-state energy $E$

\begin{rem}
\begin{itemize}
\item Our $s-$energy is twice of the $s-$energy in [KS],
$\ecal_s=2\ecal'_s$. So by the equations (\ref{Z2}) and
(\ref{K2}), we see that when $s=2$ the leading order term of the
expected average of energy is the same as the one in minimum
energy. So is the 0-energy case.

\item In equation (\ref{Z4}), we can't figure out the constant
precisely. Actually it is a conjecture in [KS]. Since in the
Green's function energy, 2-energy and 0-energy, all the leading
order terms of expected average are the same as the one in minimum
energy, this paper probably offers a method to solve the
conjecture. It will be discussed more after the proof of Theorem
1.2(1).

\end{itemize}
\end{rem}

\bigskip

 An additional motivation to study energies of random zeros is
that there are examples of numerical integration over the Riemann
surface. In numerical integration, one integrates a function with
respect to a probability measure $\mu$ by generation $N$ random
points from the ensemble $(M,\mu)$ and averaging over the points.
In this article, we generate N random points from $(M, \omega_h)$
by taking the zeros of a random polynomial. The same numerical
integration procedure is used  in the
 recent paper [DKLR] to  numerically integrate quantities  over
Calabi-Yau threefolds. The more elementary numerical integrations
in this article illustrate the speed of convergence of the
integration procedure.

\section{Background}
We begin with some notations and basic properties of sections of
holomorphic line bundles, Gaussian measures and the relation
between polynomials and sections. The notations are the same as in
[SZ1] and [BSZ]. Here we only deal with complex dimension one
case, and [PBZ] discuss the general case.

\subsection{Complex Geometry.} We denote by $(L,h)\to M$ a
holomorphic line bundle with smooth Hermitian metric $h$ whose
curvature form
\begin{equation}
\Theta_h=-\partial\bar\partial\log ||e_L||^2_h,
\end{equation}
is a positive $(1,1)$-form. Here, $e_L$ is a local non-vanishing
holomorphic section of $L$ over an open set $U\subset M$, and
$||e_L||_h=h(e_L,e_L)^{1/2}$ is the $h-$norm of $e_L$. As in
[BSZ], we give $M$ the Hermitian metric corresponding to the
K$\ddot{a}$hler form $\om=\frac{\sqrt{-1}}{2}\Theta_h$ and the
induced Riemannian volume form
\begin{equation}
dV_M=\om.
\end{equation}
We denote by $H^0(M,L^N)$ the space of holomorphic sections of
$L^N=L\otimes\cdots\otimes L.$ The metric $h$ induces Hermitian
metrics $h^N$ on $L^N$ given by $||s^{\otimes
N}||_{h^N}=||s||^N_h.$ We give $H^0(M, L^N)$ the Hermitian inner
product
\begin{equation}
<s_1,s_2>=\int_Mh^N(s_1,s_2)dV_M  \ \ \ \  (s_1,s_2\in
H^0(M,L^N)),
\end{equation}
and we write $|s|=<s,s>^{1/2}.$

For a holomorphic section $s\in H^0(M,L^N)$, we let $Z_s$ denote
the current of integration over the zero divisor of $s$:
$$(Z_s, \phi)=\int_{Z_s}\phi, \ \ \ \phi\in \dcal^{0,0}(M),$$
here $\dcal^{0,0}(\Om)$ compactly supported $(0,0)$ forms
(compactly supported smooth function) on M. A current is an
element of the dual space $\dcal^{0,0}(M)'.$

The Poincar$\acute{e}$-Lelong formula (see e.g. [GH]) expresses
the integration current of a holomorphic section $s=ge^{\otimes
N}_L$ in the form:
\begin{equation}
Z_s=\frac{i}{\pi}\partial\bar\partial\log
|g|=\frac{i}{\pi}\partial\bar\partial\log||s||_{h^N} +N\om.
\end{equation}
We also denote by $|Z_s|$ the Riemannian $0-$volume i.e.
Riemannian function along the regular points of $Z_s$, regarded as
a measure on $M$:
\begin{equation}
(|Z_s|,\phi)=\int_{Z^{reg}_s}\phi dVol;
\end{equation}

\subsection{Random sections and Gaussian measures.} We now give
$H^0(M,L^N)$ the complex Gaussian probability measure
\begin{equation}
d\mu(s)=\frac{1}{\pi^{d_N}}e^{-|c|^2}dc, \ \
s=\sum^{d_N}_{j=1}c_jS^N_j,
\end{equation}
where $\{S^N_j:1\leq j\leq d_N\}$ is an orthonormal basis for
$H^0(M,L^N)$ and $dc$ is $2d_N-$dimensional Lebesgue measure. This
Gaussian is characterized by the property that the $2d_N$ real
variable $\rcal c_j$, $\ical c_j$ ($j=1,....,d_N$) are independent
random variables with mean 0 and variance $\frac12$;i.e.,
$$\mathbf{E}c_j=0, \ \ \mathbf{E}c_jc_k=0, \ \ \mathbf{E}c_j\bar c_k=\delta_{jk}.$$
Here and throughout this article, $\mathbf{E}$ denotes
expectation: $\mathbf{E}\phi=\int\phi d\mu$.

We then regard the currents $Z_s$ (resp. measures $|Z_s|$), as
current-valued (resp. measure-valued) random variables n the
probability space ($H^0(M,L^N),d\mu$);i.e., for each test form
(resp. function) $\phi$, ($|Z_s|,\phi$) (resp. ($|Z_s|,\phi$)) is
a complex-valued random variable.

Since the zero current $Z_s$ is unchanged when $s$ is multiplied
by an element of $\C^*$, our results are the same if we instead
regard $Z_s$ as a random variable on the unit sphere $SH^0(M,L^N)$
with Haar probability measure. We prefer to use Gaussian measures
in order to facilitate computations.

\subsection{Correlation currents and measures.} The $n-$point
correlation current of the zeros is the current on $M^n=M\times
M\times...\times M$ ($n$ times) given by
\begin{equation}
K^N_n(z^1,...,z^n):=\mathbf{E}(Z_s(z^1)\otimes
Z_s(z^2)\otimes\cdots\otimes Z_s(z^n))
\end{equation}
in sense that for any test form
$\phi_1(z^1)\otimes\cdots\otimes\phi_n(z^n)\in
\dcal^{0,0}(M)\otimes\cdots\otimes\dcal^{0,0}(M).$
\begin{equation}
(K^N_n(z^1,...,z^n),\phi_1(z^1)\otimes\cdots\phi_n(z^n))=\mathbf{E}[(Z_s,\phi_1)(Z_s,\phi_2)\cdots(Z_s,\phi_n)].
\end{equation}
When $n=2$, the correlation measures take the form
\begin{equation}\label{dim 1}
 K^N_2(z,w)=[\Delta]\wedge(K^N_1(z)\otimes
1)+\kappa^N(z,w)\om_z\otimes\om_w\ \ \ (N>>0),
\end{equation}
where $[\Delta]$ denotes the current of integration along the
diagonal $\Delta={(z,z)}\subset M\times M$, and $\kappa^N\in
C^\infty(M\times M)$. In [SZ3], Bernard Shiffman and Steve
Zelditch introduced a primary object "bipotential" for the pair
correlation current; in terms of the notation used here, the
bipotential is a function $Q_N(z,w)$ such that:
\begin{equation}\label{Q}
\Delta_z\Delta_wQ_N(z,w)=K^N_2(z,w)-K^N_1(z)\wedge K^N_1(w).
\end{equation}
In [SZ3], the authors proved that for $b>\sqrt{j+2k}$, $j,k\geq
0$, we have
\begin{equation}\label{Q deri}
\nabla^jQ_N(z,w)=O(N^{-k}) \ \ \ \ uniformly \ for \ r_h(z,w)\geq
b\sqrt\frac{\log N}{N},
\end{equation}\label{2-1}
here $\nabla^jR={\frac{\partial^jR}{\partial u^{K'}\partial
u^{K"}}}: |K'|+|K"|=j$, and $r_h$ is the geodesic distance derived
by $h$. As (\ref{Q}), we have
\begin{equation}\label{K}
K^N_2(z,w)-K^N_1(z)\wedge K^N_1(w)=O(N^{-k}) \ \ \ \ uniformly \
for \ r_h(z,w)\geq b\sqrt\frac{\log N}{N}.
\end{equation}

\subsection{Relation of polynomials and sections.} By homogenizing, we may identify the space of
polynomials of degree $N$ in one complex variables with the space
$H^0(\CP^1,\ocal(N))$ of holomorphic sections of the $N-$power of
the hyperplane bundle over $\CP^1$. This space carries a natural
$SU(2)$-invariant inner product and associated Gaussian measure
$d\mu$. We associate degree $N$ polynomial $p$ zero set $Z_p=
\{p(z)=0\}$, which is almost always discrete, and thus obtain a
random point process on $\CP^1$.

\subsection{Green's function on Riemann surfaces} In this section,
we discuss Green's functions on Riemann surfaces $(M,g)$. The
Green's function is the kernel of $(-\Delta_g)^{-1}$ i.e.
$G_g(z,w)dV_g=-\Delta^{-1}$, which is orthogonal to the constant
functions, that is
\begin{equation}\label{Gp}\int_M G_g(z,w)dV_g=0
\end{equation}. Here, $-\Delta$ is the Laplacian operator.  Let $\phi_j$ be
the eigenfunctions of $-\Delta$, then
\begin{equation}
G_g(z,w)=\sum_{j\neq 0}\frac{\phi_j(z)\phi_j(w)}{\lambda^2_j},
\end{equation}
where $-\Delta \phi_j=\lambda_j\phi_j$, and $\lambda_j=0.$ So
\begin{equation}
-\Delta G_g(z,w)=\sum_{j\neq 0}
\phi_j(z)\phi_j(w)=\delta_z(w)-\frac1 {vol(M,g)}.
\end{equation}

It is well-known that $G_g(z,w)$ on Riemann surface has following
formula [see H]
\begin{equation}
G_g(z,w)\sim -\frac1{2\pi}\chi(z,w)\log r(z,w) +F(z,w)
\end{equation}
 here, $F\in C^\infty(M\times M)$ and $\chi(z,w)$ is a cut-off function which equals 1 on
 $r(z,w)\leq C_1$ and 0 on $r(z,w)\geq C_2$, where $0<C_1<C_2$.

\section{Proof of Theorem 1.1}
\begin{lem}
$\int_{M\times M}G_g(z,w)K_1^N(z)\wedge K_1^N(w)=O(N^{-2})$
\end{lem}
\begin{proof}
In [SZ3], the authors proved that
\begin{equation}
K^N_1(z)=\frac i\pi\partial\bar\partial\log \Pi_N(z,z)+\frac
N\pi\om_z,
\end{equation}
where, $\Pi_N(z,z)$ is the Szeg\"o kernel, we have following
asymptotic:
\begin{equation}
\Pi_N(z,z)=N(1+\frac {s(z)}{N}+O(N^{-2})),
\end{equation}
where $s(z)$ is the scalar curvature of $\om_z$. So we get
\begin{equation}
\log\Pi_N(z,z)=\log(N(1+\frac {s(z)}{N}+O(N^{-2}))=\log
N+\frac{s(z)}{N}+O(N^{-2})
\end{equation}
and
\begin{equation}\label{K1}
K^N_1(z)=\frac N\pi\om_z+\frac i{\pi N}\partial\bar\partial
s(z)+O(N^{-2})
\end{equation}
\begin{align*}
&\int_{M\times M}G_g(z,w)K_1^N(z)\wedge K_1^N(w)\\
=&\int_{M\times M}G_g(z,w)(\frac N\pi\om_z+\frac i{\pi
N}\partial\bar\partial s(z)+O(N^{-2}))\wedge (\frac
N\pi\om_w+\frac i{\pi
N}\partial\bar\partial s(w)+O(N^{-2}))\\
\end{align*}
Since $\int_MG_g(z,w)\om_z=0$, the last term in above equations
becomes
\begin{align*}
& \int_{M\times M}G_g(z,w)(\frac N\pi\om_z+\frac i{\pi
N}\partial\bar\partial s(z)+O(N^{-2}))\wedge (\frac
N\pi\om_w+\frac
i{\pi N}\partial\bar\partial s(w)+O(N^{-2}))\\
=&-\frac 1{(\pi N)^2}\int_{M\times M}G_g(z,w)\partial\bar\partial
s(z)\wedge \partial\bar\partial s(w)+o(N^{-2})
\end{align*}
So we have $$\int_{M\times M}G_g(z,w)K_1^N(z)\wedge
K_1^N(w)=O(N^{-2}).$$
\end{proof}
\begin{lem}
If $w=z+\frac u{\sqrt N}$, we have:
$$\E\ecal^N_{G_g} =N\int_M\int_{|u| \leq b \sqrt{\log
N}}G_g(z,z+\frac u{\sqrt N})(H(\frac12u^2)-1)\om_z\otimes\frac
i{2\pi}|u|^2 + O(N^{-k}), \;\; \forall k > 0,
$$
\end{lem}
where $H$ will be given in the proof and $|u|=\sqrt Nr_h(z,w)$.
\begin{proof}
By equation (4) and our discussion in section 2, we get:
 \begin{align*} \E\ecal^N_{G_g}&
 =\int_{H^0(M,L^N)}\ecal_{G_g}(Z_s)d\mu_h(s)\\
&=\E(G_g(z,w),Z_s \otimes Z_s - Z_{\Delta})\\
&=(G_g(z,w),K^N_2(z,w)-[\Delta]\wedge
(K^N_1(z)\otimes 1))\\
& = \int_{M \times M} G_g(z,w) (K^N_2(z,w)-[\Delta]\wedge
(K^N_1(z)\otimes 1))\\
&=\int_M\int_{r_h(z,w) \leq \frac{b \sqrt{\log N}}{\sqrt{N}}}
G_g(z,w)K^N_2(z,w)-[\Delta]\wedge (K^N_1(z)\otimes 1)\\
&+\int_M\int_{r_h(z,w) \geq \frac{b \sqrt{\log N}}{\sqrt{N}}}
G_g(z,w)K^N_2(z,w)\\
\end{align*}
By the lemma 3.1, we have:
\begin{align*}
&\int_M\int_{r_h(z,w) \leq \frac{b \sqrt{\log N}}{\sqrt{N}}}
G_g(z,w)K^N_2(z,w)-[\Delta]\wedge (K^N_1(z)\otimes 1)\\
&+\int_M\int_{r_h(z,w) \geq \frac{b \sqrt{\log N}}{\sqrt{N}}}
G_g(z,w)K^N_2(z,w)-\int_{M\times M}G(z,w)K_1^N(z)\wedge K_1^N(w)+O(N^{-2})\\
&=\int_M\int_{r_h(z,w) \leq \frac{b\sqrt{\log N}}{\sqrt{N}}}
G_g(z,w)(K^N_2(z,w)-[\Delta]\wedge
(K^N_1(z)\otimes1)-K_1^N(z)\wedge K_1^N(w) )\\
&+\int_M\int_{r_h(z,w) \geq \frac{b \sqrt{\log N}}{\sqrt{N}}}
G_g(z,w)(K^N_2(z,w)-K^N_1(z)\wedge K^N_1(w))+O(\frac 1{N^2})\\
\end{align*}

Since $G_g(z,w)= -\frac1{2\pi}\chi(z,w)\log r_g(x,y) +F(z,w)$ and
$F$ is bounded since $M$ is compact, therefore when $\
r_h(s,w)\geq b\sqrt\frac{\log N}{N}$, then $r_g(s,w)\geq
b'\sqrt\frac{\log N}{N}$, $|G_g(z,w)|$ is bounded by $\log N$, by
the  equation (\ref{2-1}),the last equation becomes:
$$\int_M\int_{r_h(z,w) \geq \frac{b \sqrt{\log N}}{\sqrt{N}}}
G_g(z,w)(K^N_2(z,w)-K^N_1(z)\wedge K^N_1(w))=O(N^{-k}),$$ so we
get
\begin{align}
\E\ecal^N_{G_g}=&\int_M\int_{r_h(z,w) \leq \frac{b \sqrt{\log
N}}{\sqrt{N}}} G_g(z,w)(K^N_2(z,w)-[\Delta]\wedge
(K^N_1(z)\otimes1)-K_1^N(z)\wedge
K^N_1(w))\notag\\
&+O(N^{-k})
\end{align}

We note there is a formula in [BSZ3] about $K_2$ on P783 Theorem
4.1

\begin{align}\label{H}
 K^N_2(z_0+\frac{z}{\sqrt N},z_0+\frac{w}{\sqrt N})\to
K^\infty_2(z,w)=&[\pi\delta_0(z-w))+H(\frac{1}{2}|z-w|^2)]\notag\\
& \cdot\frac{i}{2\pi}\partial\bar\partial
|z|^2\wedge\frac{i}{2\pi}\partial\bar\partial |w|^2
\end{align}
 where $H(t)=\frac{(sinh^2t+t^2)cosht-2tsinht}{sinh^3t}$, and
 when $t\to 0$, $H(t)=t-\frac 29t^3+O(t^5)$ and when $t\to \infty$, $H(t)=1+O(e^{-t^4})$. Here is the graph of
 $H(t)-1$ (Figure 1).
\begin{figure}[h]
\centerline{\includegraphics{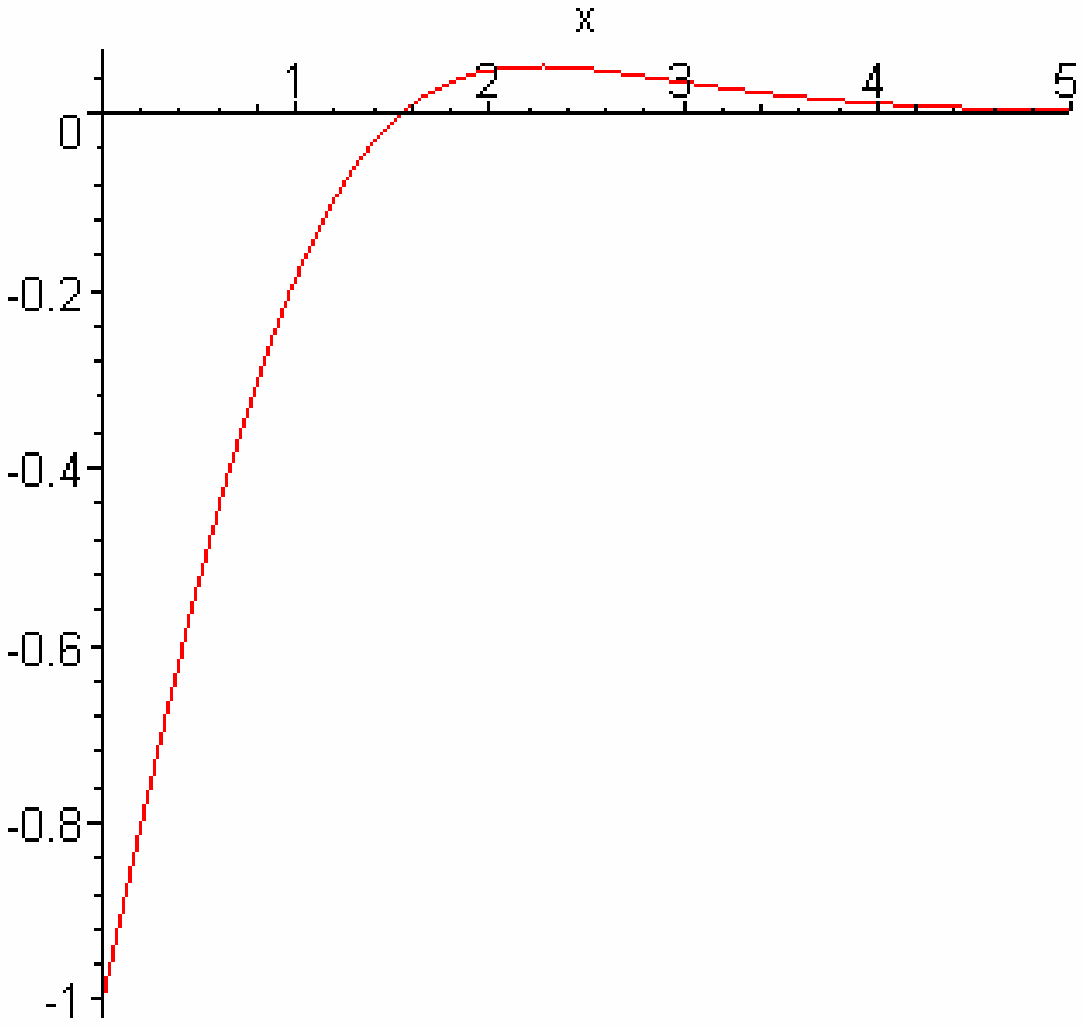}} \caption{$H(t)-1$}
\end{figure}

We change variable
$$w=z+\frac u{\sqrt N}$$
combining (\ref{dim 1}) and (\ref{H}),  we get
\begin{align}{\label{become H}}
K^N_2(z,w)-[\Delta]\wedge (K^N_1(z)\otimes1)=&K^N_2(z,z+\frac
u{\sqrt N})-[\Delta]\wedge
(K^N_1(z)\otimes1)\notag\\
=&NH(\frac12|u|^2)\cdot \frac{\om_z}\pi\wedge\frac
i{2\pi}\partial\bar\partial |u|^2+O(N^{-1})
\end{align}
Since
\begin{align}\label{kk}
K^N_1(z)\wedge K^N_1(w)&=(\frac N\pi\om_z+\frac i{\pi
N}\partial\bar\partial s(z)+O(N^{-2}))\wedge (\frac
N\pi\om_w+\frac i{\pi N}\partial\bar\partial s(w)+O(N^{-2}))\notag\\
&=\frac N\pi\om_z\wedge\frac N\pi\om_w +O(1).
\end{align}
\begin{align*} \om_w&=-\frac i2\partial\bar\partial\log ||e_L||^2_h\\
&=\frac i2\partial_w\bar\partial_w2\phi(w), \end{align*} where
$||e_L||_h=e^{-\phi(w)}$. $\phi$ is the K\"ahler potential and
since we consider the second derivative of $\phi$, so without lose
of generality we have $$\phi(z+\frac{u}{\sqrt
N})=\frac{|u|^2}{2N}+O(\frac {|u|^3}{N^{\frac 32}}).$$ Since
$u=\sqrt N(w-z)$, we have
\begin{align*}\partial_w\bar\partial_w2\phi(w)&=\partial_u\bar\partial_u2\phi(z+\frac u{\sqrt
N})+Q\\
&=\partial_u\bar\partial_u\frac{|u|^2}N+Q
\end{align*}
here, $Q$ includes $dz$ or $d\bar z$. Therefore,
\begin{align*}
\frac N\pi\om_w&=N\frac
i{2\pi}\partial_u\bar\partial_u\frac{|u|^2}N+Q\\
&=\frac i{2\pi}\partial_u\bar\partial_u|u|^2+Q
\end{align*}
 we have
\begin{equation}{\label{K-K}}
K^N_1(z)\wedge K^N_1(z+\frac u{\sqrt N})=\frac N \pi \om_z\wedge
\frac i {2\pi}\partial\bar\partial |u|^2 + O(1)
\end{equation}

 Since $r_h(z,z+\frac u{\sqrt N})=\frac{|u|}{\sqrt N}+o(N)$, we combine equations (\ref{become H}) and (\ref{K-K}) to get
\begin{align*}
\E \ecal^N_{G_g}& = \int_M \int_{r_h(z,w) \leq \frac{b \sqrt{\log
N}}{\sqrt{N}}} G_g(z,w)(K^N_2(z,w)-[\Delta]\wedge
(K^N_1(z)\otimes1)-K_1^N(z)\wedge
K^N_1(w))\\
&=N\int_M\int_{|u| \leq b \sqrt{\log N}}
G_h(z,z+\frac u{\sqrt N})(H(\frac12u^2)-1)\frac{\om_z}{\pi}\otimes\frac{i}{2\pi}\partial\bar\partial |u|^2 + O(N^{-k})\\
\end{align*}
\end{proof}

Now we complete the proof of Theorem 1.1:
\begin{proof}
According to $\S 3$
$$G_g(z,w)=-\frac 1{2\pi}\chi(z,w)\log r_g(z,w)+F(z,w).$$
$u$ is the local coordinate for $z$, therefore, if $r_h(z,z+\frac
u{\sqrt N})=\frac {|u|}{\sqrt N}+O(N^{-\frac 32})$, then
$r_g(z,z+\frac u{\sqrt N})=\frac {\sqrt {<B(z)u,u>}}{\sqrt
N}+O(N^{-\frac 32})$, where $B(z)$ is  a symmetric positive
definite operator  on $T_z M$ with respect to the metric
determined by h, once we introduce the $u$ coordinate, then $B(z)$
is a symmetric positive definite matrix which is uniformly bounded
on $M$. Then we have:
\begin{align*}
\E\ecal^N_{G_g}&=N\int_M\int_{0\leq |u|\leq b\sqrt{\log
N}}[-\frac1{2\pi}\log (\frac {\sqrt {<B(z)u,u>}}{\sqrt
N}+O(N^{-\frac 32}))][H(\frac{1}{2}|u|^2)-1]\frac{\om_z}{\pi}\otimes\frac i{2\pi}\partial\bar\partial |u|^2\\
&=\frac1{2\pi}N\log\sqrt N\int_M\frac{\om_z}{\pi}\int_{0\leq
|u|\leq b\sqrt{\log
N}}[H(\frac{1}2|u|^2)-1]\frac i{2\pi}\partial\bar\partial |u|^2\\
& \ \ \ -\frac N{2\pi}\int_M\int_{0\leq |u|\leq b\sqrt {\log
N}}\log \sqrt
{<B(z)u,u>}(H(\frac12|u|^2)-1)\frac{\om_z}\pi\otimes\frac
i{2\pi}\partial\bar\partial
 |u|^2+o(N)\\
&\sim\frac1{4\pi}N\log N\int_{0\leq |u|<
\infty}[H(\frac{1}{2}|u|^2)-1]\frac i{2\pi}\partial\bar\partial
|u|^2\\
& \ \ \ -\frac N{2\pi}\int_M\int_{0\leq |u|<\infty}\log \sqrt
{<B(z)u,u>}(H(\frac12|u|^2)-1)\frac{\om_z}\pi\otimes\frac
i{2\pi}\partial\bar\partial
|u|^2+o(N)\\
&=I+II+o(N)
\end{align*}
$\int_M\frac{\om_z}{\pi}=1$ because we assume the Chern class of
$L$, $c_1(L)=1$. Using normal coordinates, we have
\begin{align*}
&\int_{0\leq |u|<\infty}[H(\frac{1}{2}|u|^2)-1]\frac i{2\pi}\partial\bar\partial |u|^2\\
&=\frac1{\pi}\int^{2\pi}_0\int^\infty_0[H(\frac{1}{2}r^2)-1]rdrd\theta\\
&=2\int^\infty_0[H(\frac{1}{2}r^2)-1]rdr\\
&=-1
\end{align*}
 Therefore,$I=-\frac1{4\pi}N\log N.$
And since $B(z)$ is uniformly bounded on $M$, moreover, $B(z)$
varies smoothly with z  and there exist $C_1, C_2 > 0$ so that
$C_1 |u|^2\leq <B(z) u, u > \leq C_2 |u|^2$.  so it is easy to get
that $II=O(N)$.

So we have
$$\E\ecal^N_{G_g}=-\frac1{4\pi}N\log N+O(N)$$

\end{proof}
\section{Proof of Theorem 1.2}

\subsection{Proof of Theorem 1.2(1)}
\begin{proof}
\begin{align*}
\E \ecal^N_s& = \int_{\CP^1 \times \CP^1} \frac1{[z,w]^s}
(Z_s \otimes Z_s - Z_{\Delta})\\
& = \int_{\CP^1 \times \CP^1} \frac1{[z,w]^s}
(K^N_2(z,w)-[\Delta]\wedge K^N_1(z))\\
&=\int_{\CP^1}\int_{r(z,w) \leq \sqrt\frac{\log N}{N}}
[z,w]^{-s}(K^N_2(z,w)-[\Delta]\wedge
(K^N_1(z)\otimes1))\notag\\
&+\int_{\CP^1}\int_{\sqrt\frac{\log N}{N}\leq r(z,w) \leq \pi}
[z,w]^{-s}K^N_2(z,w)\notag\\
&=I+II
\end{align*}
To calculate $I$, we use the same method in $\S 3$. We change the
variables
$$w=z+\frac u{\sqrt N},$$ by equation (\ref{relat}) we get
\begin{align}\label{2.1}
I&=\int_{\CP^1}\frac{\om_z}{\pi}\int_{0\leq |u|\leq \sqrt {\log
N}}(2(1-\cos\frac {|u|}{\sqrt N}))^{-\frac s2}NH(\frac
12|u|^2)\frac
i{2\pi}\partial\bar\partial |u|^2\notag\\
&=\frac 1\pi\int^{2\pi}_0\int^{\sqrt{\log N}}_0(2(1-\cos\frac
{r}{\sqrt N}))^{-\frac s2}NH(\frac 12r^2)rdrd\theta
\end{align}
Since $2(1-\cos\frac r{\sqrt N})\sim \frac {r^2}N$, then
(\ref{2.1}) becomes
\begin{align}\label{2.1.1}
&2N^{1+\frac s2}\int^{\sqrt{\log N}}_0\frac 1{r^s}H(\frac12
r^2)rdr\notag\\
=&2N^{1+\frac s2}\int^{\sqrt{\log N}}_{\frac ML}H(\frac12
r^2)r^{1-s}dr+2N^{1+\frac s2}\int^{\frac M L}_0H(\frac12
r^2)r^{1-s}dr
\end{align}
Here,$\frac ML<\sqrt{\log N}$ then we can assume $\frac
ML=O(\sqrt[4]{\log N})$, since $H(\frac12r^2)\to \frac 12r^2$ as
$r\to 0$,then let $L\to \infty$,
\begin{itemize}
\item When $s=2$, the second part of (\ref{2.1.1}) is asymptotic
to \begin{equation} \frac 12N^2(\frac ML)^2.
\end{equation}

\item When $s<4$ and $s\neq 2$, the second part of (\ref{2.1.1})
is asympotitc to
\begin{equation}
\frac {N^{1+\frac s2}}{4-s}(\frac ML)^{4-s}.
\end{equation}
\end{itemize}

 And since $H(r)\to 1$ as $r\to\infty$, then let
$M\to\infty$,
\begin{itemize}
\item When $s=2$, the first part of (\ref{2.1.1}) is asymptotic to
\begin{equation}
N^2\log(\log N)-2N^2\log(M/L).
\end{equation}

\item When $s<4$ and $s\neq 2$, the first part of (\ref{2.1.1}) is
asymptotic to \begin{equation}\frac{N^{1+\frac s2}}{2-s}(\log
N)^{1-\frac s2}+N^{1+\frac s2}(\frac ML)^{2-s}.
\end{equation}
\end{itemize}

So \begin{itemize} \item When $s=2$,
\begin{equation}\label{Is=2}
I=N^2\log(\log N)+\frac 12N^2(\frac ML)^2-2N^2\log(M/L).
\end{equation}

\item When $s<4$ and $s\neq 2$,
\begin{equation}\label{Is<4}
I=\frac {N^{1+\frac s2}}{4-s}(\frac ML)^{4-s}+\frac{N^{1+\frac
s2}}{2-s}(\log N)^{1-\frac s2}+N^{1+\frac s2}(\frac ML)^{2-s}.
\end{equation}
\end{itemize}

To calculate $II$, we use the equation (\ref{K})
\begin{equation*}
K^N_2(z,w)-K^N_1(z)\wedge K^N_1(w)=O(N^{-k}) \ \ \ \ uniformly \
for \ r_h(z,w)\geq b\sqrt\frac{\log N}{N}.
\end{equation*}
and equation (\ref{kk})
\begin{equation*}K^N_1(z)\wedge K^N_1(w)=\frac N \pi \om_z\wedge \frac N \pi \om_w+ O(1)
\end{equation*}
to get \begin{align}\label{II}
II&=\int_{\CP^1}\int_{\sqrt\frac{\log N}{N}\leq r(z,w)\leq
2}[z,w]^{-s}
K^N_1(z)\wedge K^N_1(w)+O(N^{-k})\notag\\
&=N^2\int_{\CP^1}\frac{\om_z}{\pi}\int_{\sqrt\frac{\log N}{N}\leq
r(z,w)\leq 2}[z,w]^{-s}\frac {\om_w}\pi+O(1)
\end{align}
Since $\int_{\CP^1}\frac\om\pi=1$,  if we use azimuthal angle
$\phi$, we get $\int_{S^2}\sin\phi d\phi d\theta=4\pi$.  For the
standard unit sphere, $\phi=r$, where $r$ is the round distance.
Then we have:
\begin{align}
II&=\frac
{N^2}{4\pi}\int^{2\pi}_0\int^\pi_{\sqrt\frac{\log N}{N}}\frac 1{(2(1-\cos\phi))^{\frac s2}}\sin\phi d\phi d\theta\notag\\
&=\frac {N^2}{2}\int^\pi_{\sqrt\frac{\log N}{N}}\frac
1{(2(1-\cos\phi))^{\frac s2}}\sin\phi d\phi\notag\\
&=\frac{N^2}{4}\int^\pi_{\sqrt\frac{\log N}{N}}\frac
1{(2(1-\cos\phi))^{\frac s2}}d(2(1-\cos\phi))
\end{align}

\begin{itemize}
\item When $s=2$, \begin{align}\label{IIs=2} II=& \frac
{N^2}4\log(2(1-\cos\phi))|^\pi_{\sqrt\frac{\log N}{N}}\notag\\
=& \frac {N^2}2\log 2-\frac {N^2}4\log(2(1-\cos\sqrt\frac{\log
N}{N}))\notag\\
\sim&\frac {N^2}4\log N-\frac {N^2}4\log(\log N)+\frac {N^2}2\log
2.
\end{align}

\item When $s<4$ and $s\neq 2$
\begin{align}\label{IIs<4}II=&\frac2{2-s}\frac{N^2}4(2(1-\cos\phi))^{1-\frac
s2}|^\pi_{\sqrt\frac{\log N}{N}}\notag\\
=&\frac{2^{1-s}N^2}{2-s}-\frac{N^2}{2(2-s)}(2(1-\cos\sqrt\frac{\log
N}{N}))^{1-\frac s2}\notag\\
\sim&\frac{2^{1-s}N^2}{2-s}-\frac {N^2}{2(2-s)}(\frac{\log
N}{N})^{1-\frac s2}\notag\\
=&\frac{2^{1-s}N^2}{2-s}-\frac1{2(2-s)}N^{1+\frac s2}(\log
N)^{1-\frac s2}
\end{align},
\end{itemize}

\begin{itemize}
\item When $s=2$, since $\frac ML<\sqrt{\log N}$, then
$$\E\ecal^N_2=\frac14{N^2}\log N+\frac {3N^2}4\log(\log N)+\frac {N^2}2\log
2+\frac12N^2(\frac ML)^2-2N^2\log\frac ML+o(N^2).$$

\item When $s<2$,
$$\E\ecal^N_{s<2}=\frac{2^{1-s}}{2-s}N^2+\frac 1{(2(2-s))}N^{1+\frac s2}(\log N)^{1-\frac s2}+o(N^{1+\frac s2}(\log N)^{1-\frac s2})$$

\item When $2<s<4$, the leading order term is $\frac {N^{1+\frac
s2}}{4-s}(\frac ML)^{4-s}$ in (\ref{Is<4}), however it is hard to
figure out what $\frac ML$ is.
$$\E\ecal^N_{2<s<4}=C\frac {N^{1+\frac s2}}{4-s}+O(N^{1+\frac s2}(\log N)^{1-\frac s2}).$$
\end{itemize}
\end{proof}
When $2<s<4$, it is hard for us to figure out the constant $C$,
because we can't give the asymptotic to the integration in
(\ref{2.1}).

\subsection{Proof of Theorem 1.2(2)}

\begin{proof}
\begin{align}\label{log}
\E \ecal^N_0& = \int_{\CP^1 \times \CP^1}
-\log [z,w](Z_s \otimes Z_s - Z_{\Delta})\notag\\
& = \int_{\CP^1 \times \CP^1} -\log [z,w]
(K^N_2(z,w)-[\Delta]\wedge K^N_1(z))\notag\\
&=\int_{\CP^1}\int_{r(z,w)\leq \frac{\sqrt{\log N}}{\sqrt{N}}}
-\log [z,w](K^N_2(z,w)-[\Delta]\wedge
(K^N_1(z)\otimes1))\notag\\
&+\int_{\CP^1}\int_{\frac{\sqrt{\log N}}{\sqrt{N}}\leq r(z,w) \leq
\pi}
-\log [z,w]K^N_2(z,w)\notag\\
&=I+II
\end{align}
As \S 4.1, we change variables
$$w=z+\frac u{\sqrt N}$$ and by equation (\ref{relat}), we get
\begin{align}\label{2.2}
I&=\int_{\CP^1}\frac{\om_z}{\pi}\int_{0\leq |u|\leq \sqrt {\log
N}}-\log\sqrt{2(1-\cos\frac{|u|}{\sqrt N})}NH(\frac 12|u|^2)\frac
i{2\pi}\partial\bar\partial |u|^2\notag\\
&=\frac 1\pi\int^{2\pi}_0\int^{\sqrt{\log N}}_0-\log\sqrt{2(1-\cos\frac{r}{\sqrt N})}NH(\frac 12r^2)rdrd\theta\notag\\
&=-N\int^{\sqrt{\log N}}_0(\log2)H(\frac
12r^2)rdr-N\int^{\sqrt{\log
N}}_0\log(1-\cos\frac{r}{\sqrt N})H(\frac 12r^2)rdr\\
\end{align}
Since $H(r)\to r$ as $r\to 0$ and $H(r)\to 1$ as $r\to \infty$, we
get:

$$I\sim-\frac{\log2}2N\log N-N\int^{\sqrt{\log
N}}_0\log(1-\cos\frac{r}{\sqrt N})H(\frac 12r^2)rdr.$$ Since
$1-\cos\frac r{\sqrt N}=\frac {r^2}{2N}$
\begin{align*}
&\int^{\sqrt{\log N}}_0\log(1-\cos\frac{r}{\sqrt
N})H(\frac 12r^2)rdr\\
=&\int^{\sqrt{\log N}}_0\log(\frac{r^2}{2N})H(\frac 12r^2)rdr\\
=&\int^{\sqrt{\log N}}_0(\log\frac{r^2}2)H(\frac
12r^2)d\frac{r^2}2-\int^{\sqrt{\log
N}}_0\log(N)H(\frac 12r^2)rdr\\
=&\frac12\log(\log N)\log N-\frac12(\log2+1)\log N-\frac12(\log2+1)-\frac{\log^2N}2\\
\end{align*}
So $I=\frac N2\log^2N-\frac12N\log(\log N)\log N+\frac12N\log
N+\frac12(\log2+1)N$

  To calculate $II$, we use the same method in \S4.1 and get
 \begin{align}\label{II}
II&=\int_{\CP^1}\int_{\sqrt\frac{\log N}{N}\leq r(z,w)\leq
\pi}-\log
[z,w]K^N_1(z)\wedge K^N_1(w)+O(N^{-k})\notag\\
&=N^2\int_{\CP^1}\frac{\om_z}{\pi}\int_{\sqrt\frac{\log N}{N}\leq
r(z,w)\leq \pi}-\log [z,w]\frac {\om_w}\pi+O(1).
\end{align}
Since $\int_{\CP^1}\frac\om\pi=1$,  if we use azimuthal angle
$\phi$, we get $\int_{S^2}\sin\phi d\phi d\theta=4\pi$.  For the
standard unit sphere, $\phi=r$, where $r$ is the round distance.
Then (\ref{II}) becomes
\begin{align}\label{II'}
&\frac {N^2}
{4\pi}\int^{2\pi}_0\int_{\sqrt\frac{\log N}{N}\leq \phi\leq \pi}-\log\sqrt{2(1-\cos\phi)}\sin\phi d\phi d\theta\\
=&-\frac {N^2}{4}\int^\pi_{\sqrt {\frac{\log
N}{N}}}\log2(1-\cos\phi)\sin\phi d\phi\notag\\
=&-(\log2-\frac12)N^2+O(\log N)
\end{align}
In the end, we get \begin{align*}\E
\ecal^N_0=&-(\log2-\frac12)N^2+\frac N2\log^2N-\frac12N\log(\log
N)\log N\\
&+\frac12N\log N+\frac12(\log2+1)N+o(N).\\
\end{align*}
\end{proof}

\begin{appendix}
\section{} In the appendix, we give a picture
which describes the distribution to random zeros of a given random
polynomial. Let
$$p(z)=\sum^{50}_{i=1}c_i(C^j_N)^\frac12z^i,$$ where $\E(c_i)=0$ and
$\E(|c_i|^2)=1.$
\begin{figure}[h]
\centerline{\includegraphics{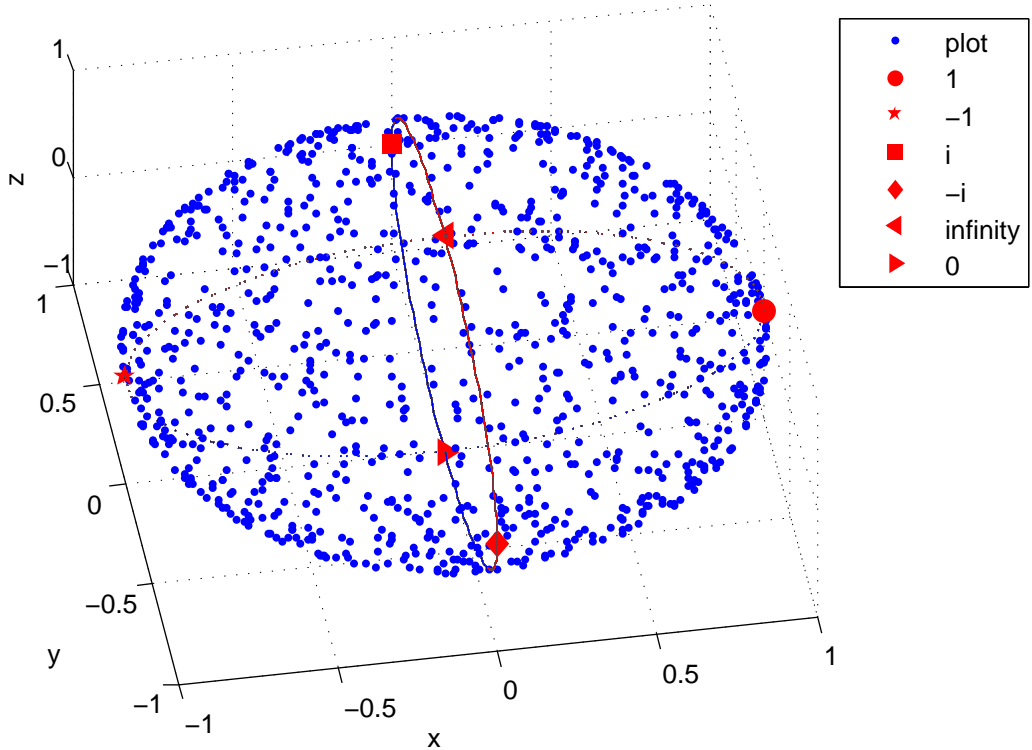}} \caption{Distribution of
zeros}
\end{figure}
\end{appendix}

\end{document}